\documentclass[twoside,leqno]{article}

\usepackage[letterpaper]{geometry}

\usepackage{ltexpprt}
\usepackage{hyperref}
 \usepackage{float}  
\usepackage{subcaption}
\usepackage{wrapfig} 
\usepackage{nicefrac} 
\usepackage{subcaption}
 
 \usepackage{amsmath}
 \usepackage{amssymb}
\usepackage{ulem}
\usepackage{natbib}
\usepackage{hyperref}
\hypersetup{
  colorlinks = true,
  citecolor  = blue,
  linkcolor  = red,
  filecolor  = blue,
  urlcolor   = blue,
}
\usepackage{lipsum}
\usepackage{amsfonts}
\usepackage{graphicx}
\usepackage{epstopdf}
\usepackage{algorithmic}

\usepackage{tikz} 
\usetikzlibrary{calc,intersections,through,arrows,angles, quotes}
\usepackage{tkz-euclide} 
\usepackage{nicefrac} 
\usepackage{subcaption}
\usepackage{wrapfig} 
\usepackage{mdframed} 

\usepackage{array}  

\ifpdf
  \DeclareGraphicsExtensions{.eps,.pdf,.png,.jpg}
\else
  \DeclareGraphicsExtensions{.eps}
\fi

\newcommand{\inp}[2]{\left\langle#1, #2 \right\rangle}
\newcommand{\norm}[1]{\left\lVert#1\right\rVert}

\newcommand{\OO}[1]{O\left( #1\right)}
\newcommand{\re}{\mathbb{R}}

\newcommand{\ee}{\eta}
\newcommand{\eep}{\widetilde{\eta}}

\newcommand{\xs}{x_*} 
\newcommand{\eeo}{\text{($\mathcal{E}^{\mathsf{GD}}_1$)}}
\newcommand{\eet}{\text{($\mathcal{E}^{\mathsf{GD}}_2$)}}
\newcommand{\ffo}{\text{($\mathcal{E}^{\mathsf{AGM}}_1$)}}
\newcommand{\fft}{\text{($\mathcal{E}^{\mathsf{AGM}}_2$)}}
\newcommand{\ggo}{\text{($\mathcal{E}^{\mathsf{SIM'}}_1$)}}
\newcommand{\ggt}{\text{($\mathcal{E}^{\mathsf{SIM'}}_2$)}}
\newcommand{\hho}{\text{($\mathcal{E}^{\mathsf{SIM}}_1$)}}
\newcommand{\hht}{\text{($\mathcal{E}^{\mathsf{SIM}}_2$)}}
\newcommand{\fsi}{f^{\Psi}}

\newcommand{\pss}{\phi^{\Psi}}
\newcommand{\low}[2]{\mathsf{LOWER}({#1};{#2})}
\newcommand{\upp}[2]{\mathsf{UPPER}({#1};{#2})}

\DeclareMathOperator*{\argmin}{argmin}

\newcommand{\breg}[3]{D_{#1}\left(#2,#3\right)}

\newcommand{\cca}[2]{\alpha^{(#1)}_{#2}}
\newcommand{\ccb}[2]{\beta^{(#1)}_{#2}}
\newcommand{\ccc}[2]{\gamma^{(#1)}_{#2}}
\newcommand{\ccd}[2]{\delta^{(#1)}_{#2}}

\begin{document}

%



\title{\Large  Understanding Nesterov's Acceleration via Proximal Point Method}
\author{Kwangjun Ahn\thanks{Department of Electrical Engineering and Computer Science,  Massachusetts Institute of Technology,  \texttt{kjahn@mit.edu}. This work was done as the author's class project for 6.881  Optimization for Machine Learning at MIT, Spring 2020.} \and Suvrit Sra\thanks{Department of Electrical Engineering and Computer Science,  Massachusetts Institute of Technology,  \texttt{suvrit@mit.edu}.} }

\date{}

\maketitle


\fancyfoot[R]{\scriptsize{Copyright \textcopyright\ 2022 by SIAM\\
Unauthorized reproduction of this article is prohibited}}





\begin{abstract} \small\baselineskip=10pt
The proximal point method (PPM) is a fundamental  method in optimization that is often used as a building block for designing optimization algorithms. In this work, we use the PPM method to  provide conceptually simple derivations along with convergence analyses of different versions of Nesterov's accelerated gradient method (AGM). The key observation is that AGM is a simple approximation of PPM, which results in an elementary derivation of the update equations and stepsizes of AGM. This view also leads to a transparent and conceptually simple analysis of AGM's convergence by using the analysis of PPM. The derivations also naturally extend to the strongly convex case. Ultimately, the results presented in this paper are of both didactic and  conceptual value; they unify and explain existing variants of AGM while motivating other accelerated methods for practically relevant settings.
\end{abstract}

\section{Introduction}
 In 1983, Nesterov introduced the \emph{accelerated gradient method} (AGM) for minimizing a convex function $f:\re^d \to \re$~\citep{nesterov1983method}.
The remarkable property of AGM is that AGM  achieves a strictly faster convergence rate  than the standard gradient descent (GD). Assuming that $f$ has Lipschitz continuous gradients, $T$ iterations of AGM are guaranteed to output a point $x_T$ with the suboptimality gap $f(x_T)-\min_x f(x) \leq O(1/T^2)$, whereas GD only ensures a suboptimality gap of $O(1/T)$. On top of being a landmark result of convex optimization, AGM is easy to implement and has found value in a myriad of applications such as sparse linear regression~\citep{beck2009fast}, compressed sensing~\citep{becker2011nesta}, the maximum  flow problem~\citep{lee2013new}, and deep neural networks~\citep{sutskever2013importance}.

AGM's importance to both theory and practice has led to a flurry of works that seek to understand its scope and the principles that underlie it~\citep{su2014differential,krichene2015accelerated,wibisono2016variational,lessard2016analysis, wilson2016lyapunov, allen2014linear, diakonikolas2019approximate}; see \S\ref{sec:related} for details. However, one curious aspect of AGM that is not yet well-understood is the fact that it appears in various different forms. Below, we list the four most representative ones:
\vspace{-10pt}
 \begin{center}
     \begin{minipage}{0.45\textwidth} 
      \begin{align*} 
    \begin{split}
 & \textstyle z_{t+1} = y_{t}-\cca{1}{t}\nabla f(y_{t})\,, \\
& \textstyle y_{t+1} = z_{t+1} +\ccb{1}{t}(z_{t+1}-z_t)\,. \end{split}
  \end{align*} 
  \begin{center}
Form I~\citep{nesterov1983method,beck2009fast}.      
  \end{center}

\end{minipage}
 \begin{minipage}{0.45\textwidth} 
    \begin{align*}
    \begin{split}
 & \textstyle y_{t} = \cca{2}{t} x_t+(1-\cca{2}{t}) z_t \,,  \\
& \textstyle z_{t+1} = y_t -\ccb{2}{t} \nabla f(y_{t})\,, \\
& \textstyle x_{t+1} = x_{t}-\ccc{2}{t} \nabla f(y_{t})\,. 
\end{split}
  \end{align*}   
    \begin{center}
Form II~\citep{nesterov2018lectures,allen2014linear}.      
  \end{center}
\end{minipage}
 \end{center}
 \begin{center}
 \begin{minipage}{0.45\textwidth} 
    \begin{align*} 
    \begin{split}
  & \textstyle y_{t} = \cca{3}{t} x_t+ (1-\cca{3}{t}) z_t\,,  \\
& \textstyle x_{t+1} = x_t -\ccb{3}{t} \nabla f(y_{t})\,, \\
 & \textstyle z_{t+1} = \ccc{3}{t} x_{t+1}+ (1-\ccc{3}{t}) z_t\,.  
\end{split}
  \end{align*}   
      \begin{center}
Form III~\citep{auslender2006interior,tseng2008accelerated,gasnikov2018universal}.      
  \end{center}
\end{minipage}
  \begin{minipage}{0.45\textwidth} 
       \begin{align*}  
   & \textstyle y_{t} =\cca{4}{t} x_t + (1-\cca{4}{t}) z_t\,, \\ 
& \textstyle x_{t+1} =  \ccb{4}{t} x_t + (1-\ccb{4}{t})y_t-  \ccc{4}{t} \nabla f(y_{t})\,, \\
& \textstyle z_{t+1} = y_{t}-\ccd{4}{t} \nabla f(y_{t})\,.          
   \end{align*} 
      \begin{center}
Form IV~\citep{nesterov2018lectures}.      
  \end{center}
\end{minipage} 
 \end{center}
The parameters $\cca{\cdot}{t},\ccb{\cdot}{t},\ccc{\cdot}{t}, \ccd{\cdot}{t}$ are stepsizes that are carefully chosen to ensure an accelerated rate. An immediate question that one may ask is: \emph{can we understand these variants of AGM in a unified manner?}

This paper answers this question by developing a transparent and unified analysis that captures all these variants of AGM  by connecting them to the proximal point method (PPM). PPM is a well-known optimization algorithm that is often used as a conceptual building block for designing other optimization algorithms  (see \S\ref{sec:ppm} for more background). The key insight (presented in \S\ref{sec:deriv}) is that one can obtain AGM simply by viewing it as an approximation of PPM. This insight is inspired by the approach of \citet{defazio2019curved}, but now with more general acceleration settings and importantly, without any recourse to duality.  
  
\paragraph{Contributions.} In summary, we make the following contributions:
\begin{itemize}
  \setlength{\itemsep}{0pt}
    \item We present an intuitive derivation of AGM by viewing it as an approximation of the proximal point method (PPM), a foundational, classical method in optimization.
    \item We present a unified method for deriving different versions of AGM, which may be of wider pedagogical interest. In particular, our approach readily extends to the strongly convex case and offers a short derivation of the most general version of AGM introduced by Nesterov in his textbook~\citep[(2.2.7)]{nesterov2018lectures}.
\end{itemize}
  
We believe that the simple derivations presented in this paper are not only of pedagogical value but are also helpful for research because they clarify, unify, and deepen our understanding of the phenomenon of acceleration. The PPM view offers a transparent analysis of AGM based on the convergence analysis of PPM~\citep{guler1991convergence}. Moreover, as we present in \S\ref{sec:tri}, the PPM view also motivates the key idea of the \emph{method of similar triangles}, a version of AGM shown to have important extensions to practically relevant settings~\citep{tseng2008accelerated,gasnikov2018universal}. Our approach also readily extends to the strongly convex case (\S\ref{sec:str}). Finally, since PPM has been studied in settings much wider than convex optimization (see e.g.,~\citep{bacak2014convex}), we believe the connections exposed herein will help in advancing the development of accelerated methods in those settings.

\noindent Before presenting our derivations, let us first recall a brief background on the proximal point method.

\section{Brief background on the proximal point method}\label{sec:ppm}
The \emph{proximal point method (PPM)}~\citep{moreau1965proximite,martinet1970regularisation,rockafellar1976monotone} is a fundamental method in optimization which solves the minimization of the cost function $f:\re^d\to \re$ by iteratively solving the  subproblem
\begin{align} 
\label{prox}
  x_{t+1} \leftarrow \argmin_{x\in \re^d} \left\{ f(x) + \frac{1}{2\ee_{t+1}} \norm{x-x_t}^2 \right\}
\end{align}
for a stepsize $\ee_{t+1}>0$, where the norm is chosen as the $\ell_2$ norm.
Despite its simplicity, solving \eqref{prox} is in general as difficult as solving the original optimization problem, and PPM is largely regarded as a ``conceptual'' guiding principle for accelerating optimization algorithms~\citep{drusvyatskiy2017proximal}. 

The baseline of our discussion is the following convergence rate of PPM for convex costs  proved  in a seminal paper by \citet{guler1991convergence}  (here $\xs$ denotes a global optimum point, i.e., $\xs \in \argmin_x f(x)$):
\begin{align}
    \boxed{\textstyle f(x_{T})- f(\xs) \leq \OO{\big(\sum_{t=1}^{T} \ee_t\big)^{-1}}\quad \text{for any $T\geq 1$.}} \label{conv:prox} 
\end{align}
In words, one can achieve an arbitrarily fast convergence rate by choosing stepsizes $\ee_t$'s large.
Below, we review  a short Lyapunov function proof of \eqref{conv:prox}, which will serve as a backbone to other analyses.
\begin{proof}[{\bf Proof of \eqref{conv:prox}}]
It turns out that the following Lyapunov function is suitable:
\begin{align}
    \boxed{\textstyle \Phi_t:=  \big(\sum_{i=1}^{t} \ee_i \big)\cdot \big(f(x_t)-f(\xs)\big) + \frac{1}{2}\norm{\xs-x_t}^2,}\label{def:lya}
\end{align}
where $\Phi_0:= \frac{1}{2}\norm{\xs-x_0}^2$ and here and below, $\norm{\cdot}$ is the $\ell_2$ norm unless stated otherwise.
Now,  it suffices to show that  $\Phi_t$  is decreasing, i.e., $\Phi_{t+1}\leq \Phi_t$ for all $t\geq 0$.
Indeed, if $\Phi_t$ is decreasing, we have $\Phi_T\leq \Phi_0$ for any $T\geq 1$, which  precisely recovers \eqref{conv:prox}. To that end, we use a standard result: 
\begin{proposition}[Proximal inequality~(see {\citep[Proposition 12.26]{bauschke2011convex}})] \label{prop:per} For a convex function $\phi:\re^d\to \re$,  let $x_{t+1}$ be the unique minimizer of the following proximal step:  $x_{t+1} \leftarrow \argmin_{x\in \re^d} \left\{\phi(x) +\frac{1}{2}\norm{x-x_t}^2\right\}$. Then, for any $u\in \re^d$, \begin{align*}
    \phi(x_{t+1})-\phi(u) +\frac{1}{2}\norm{u-x_{t+1}}^2 +\frac{1}{2}\norm{ x_{t+1}-x_t}^2  -\frac{1}{2}\norm{u-x_t}^2\leq 0\,.
\end{align*}
\end{proposition}    
Now Proposition~\ref{prop:per} completes the proof as follows:  First, we apply Proposition~\ref{prop:per} with $\phi=\ee_{t+1} f$ and $u=\xs$ and drop the term $\frac{1}{2}\norm{x_{t+1}-x_t}^2$ to obtain:
    \begin{align}
    \tag{\text{$\mathsf{Ineq}_1$}}  \ee_{t+1}\left[f( x_{t+1}) -f( \xs)\right] + \frac{1}{2} \norm{\xs-    x_{t+1}}^2  -\frac{1}{2} \norm{ \xs-x_t}^2\leq  0\,.\label{ineq:1}
\end{align} 
Next, from the optimality of $x_{t+1}$, it readily follows that  
    \begin{align}
    \tag{\text{$\mathsf{Ineq}_2$}} f(x_{t+1}) - f(x_t) \leq 0\,. \label{ineq:2}
    \end{align}
Now,  computing $\eqref{ineq:1}+( \sum_{i=1}^{t} \ee_i )\times $\eqref{ineq:2} yields $\Phi_{t+1}\leq \Phi_t$, which finishes the proof. \end{proof}

\subsection{Our conceptual question}
Although the convergence rate \eqref{conv:prox} seems powerful, it does not have any practical values as PPM is in general not implementable. 
Nevertheless, one can ask the following conceptual question:
\begin{center}
    \emph{``Can we efficiently approximate PPM for a large stepsize $\ee_t$?''}
\end{center} 
Perhaps, the most straightforward approximation would be to replace the cost function $f$ in \eqref{prox} with its lower-order approximations.
We implement this idea in the next section.

\section{Two simple approximations of the proximal point method} \label{sec:warmup}

To analyze approximation errors, let us assume that the cost function $f$ is $L$-smooth.
\begin{definition}[Smoothness] \label{def:ell2}
For $L>0$, we say a differentiable function $f:\re^d\to \re$ is  $L$-smooth if  $f(x) \leq f(y) + \inp{\nabla f(y)}{x-y} +\frac{L}{2}\norm{x-y}^2 $ for any $x,y\in\re^d$.
\end{definition}
\noindent From the convexity and the $L$-smoothness of $f$, we have the following lower and upper  bounds: for any $x,y\in \re^d$,
\begin{align*}
 \boxed{\underbrace{f(y) +\inp{\nabla f(y)}{x-y}}_{{    \textstyle=:\low{x}{y}}}    \leq f(x) \leq \underbrace{f(y) +\inp{\nabla f(y)}{x-y}+\frac{L}{2}\norm{x-y}^2}_{\textstyle =:\upp{x}{y}}\,.}
\end{align*}
In this section, we use these bounds to approximate PPM. 
\subsection{First approach: using first-order approximation} 
\label{sec:app1}
\begin{figure}
\begin{subfigure}{.49\textwidth}
  \centering
  \includegraphics[width=\linewidth]{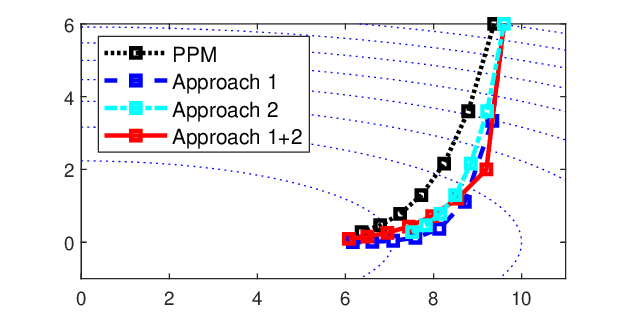} \caption{$\ee_t\equiv 1/3$.}
\end{subfigure}
\begin{subfigure}{.49\textwidth}
  \centering
  \includegraphics[width=\linewidth]{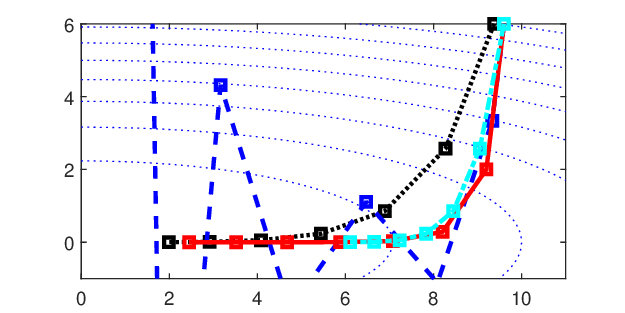}  \caption{$\ee_t=t/3$.}
\end{subfigure}
\caption{ Iterates comparison between PPM~\eqref{prox}, the first approach~\eqref{gd}, the second approach~\eqref{gd:2-1}, and the combined approach~\eqref{approx:ppm}. For the setting, we choose $f(x,y) = 0.1x^2+y^2$ and   $x_0=(10,10)$. }
\label{fig:1} 
\end{figure}

Let us first replace $f$ in the objective \eqref{prox} with its lower approximation: 
\begin{align}
   x_{t+1} \leftarrow \argmin_{x} \left\{  \low{x}{x_t}  + \frac{1}{2\ee_{t+1}} \norm{x-x_t}^2 \right\}\,. \label{gd}
\end{align}
Writing the optimality condition, one quickly notices that \eqref{gd} actually leads to  gradient descent: 
\begin{align}
x_{t+1} =x_t -\ee_{t+1} \nabla f(x_{t})\,.\label{gd:1}
\end{align}
Let us see how well \eqref{gd} approximates PPM:

\begin{proof}[{\bf Analysis of the first approach}] We first establish counterparts  of  \eqref{ineq:1} and \eqref{ineq:2}.
First, we apply Proposition~\ref{prop:per} with $\phi(x)=\ee_{t+1}\low{x}{x_t}$ and $u=\xs$:
    \begin{align*}
        &\phi(x_{t+1})-\phi(\xs)+ \frac{1}{2} \norm{\xs -x_{t+1}}^2 + \frac{1}{2} \norm{ x_{t+1}-x_t}^2 -\frac{1}{2 } \norm{\xs -x_{t}}^2\leq 0\,.
    \end{align*} 
    Now using convexity and $L$-smoothness, we have 
    \begin{align*}
        \phi(x)\leq \ee_{t+1} f(x)\leq \phi(x) +\frac{L\ee_{t+1}}{2}\norm{x-x_t}^2\,,
    \end{align*} and hence the above inequality implies the following analogue of \eqref{ineq:1}:
    \begin{align}  \tag{\text{$\mathsf{Ineq}^{\mathsf{GD}}_1$}} \ee_{t+1} \left[ f(x_{t+1})- f(\xs)\right] + \frac{1}{2} \norm{\xs -x_{t+1}}^2- \frac{1}{2 } \norm{\xs -x_{t}}^2 \leq   \eeo,  \label{ineq:1:gd}
    \end{align}
    where $\eeo:= (\frac{L\ee_{t+1}}{2}-\frac{1  }{2} )\norm{x_{t+1}-x_t}^2$.
Next, we use the $L$-smoothness of $f$ and the fact $\nabla f(x_t)= -\nicefrac{1}{\ee_{t+1}}(x_{t+1}-x_t)$ (due to \eqref{gd:1}), to obtain the following analogue of \eqref{ineq:2}:
    \begin{align}  
 \tag{\text{$\mathsf{Ineq}^{\mathsf{GD}}_2$}}    f(x_{t+1}) -f(x_t) \leq  \inp{\nabla f(x_t)}{x_{t+1}-x_t} + \frac{L}{2} \norm{x_{t+1}-x_t}^2 =   \eet, \label{ineq:2:gd}
    \end{align}
    where $\eet:=( \frac{L}{2}- \frac{1}{\ee_{t+1}}) \norm{x_{t+1}-x_t}^2$.

Now paralleling the proof of \eqref{conv:prox},  to show that $\Phi_t$~\eqref{def:lya} is a valid Lyapunov function, we need to find the stepsizes $\ee_t$'s that satisfy the following relation: $\eeo + ( \sum_{i=1}^t \ee_i) \times \eet \leq 0$. 
On the other hand, note that both \eeo~and \eet~become positive numbers when $\ee_{t+1}>2/L$.
Hence, the admissible choices for $\ee_t$ at each iteration are upper bounded by $2/L$, which together with the PPM convergence rate \eqref{conv:prox} implies that $O(\nicefrac{1}{\sum_{t=1}^T\ee_t}) =O(\nicefrac{1}{T})$ is the best convergence rate  one can prove.
Indeed, choosing  $\ee_t\equiv 1/L$, then we have $\eeo=0$ and $\eet<0$, obtaining the well-known bound of $f(x_T)-f(\xs) \leq \frac{L\norm{x_0-\xs}^2}{2T} =O(\nicefrac{1}{T})$. 
\end{proof} 
\noindent To summarize, the first approach only leads to a disappointing result: the approximation is valid only for the small stepsize regime of $\ee_t = \OO{1/L}$. 
We empirically verify this fact for a  quadratic cost in Figure~\ref{fig:1}.
As one can see from Figure~\ref{fig:1}, the lower approximation approach \eqref{gd}  overshoots for large stepsizes like  $\ee_t = \Theta(t)$ and quickly steers away from  the PPM iterates.

\subsection{Second approach: using smoothness} \label{sec:app2}
After seeing the disappointing outcome of the first approach, our second approach is to replace $f$ with  its upper approximation  due to the $L$-smoothness:
\begin{align}
     x_{t+1} \leftarrow \argmin_{x} \left\{ \upp{x}{x_t}+ \frac{1}{2\ee_{t+1}} \norm{x-x_t}^2 \right\}\,. \label{gd:2-1}
\end{align}
Writing the optimality condition,  \eqref{gd:2-1} actually leads to a conservative update of gradient descent: 
\begin{align}
x_{t+1} = x_t -\frac{1}{L+\ee_{t+1}^{-1}} \nabla f(x_{t})\,. \label{gd:2-2}
\end{align}
Note that regardless of how large $\ee_{t+1}$ we choose, the actual update stepsize in \eqref{gd:2-2} is always upper bounded by $\nicefrac{1}{L}$.
Although this conservative update  prevents the  overshooting phenomenon of the first approach, as we increase $\ee_{t}$, this conservative update becomes too tardy to be a good approximation of PPM; see Figure~\ref{fig:1}.

\section{Nesterov's acceleration via alternating two approaches}
\label{sec:deriv}
In the previous section, we have seen that the two simple approximations of PPM both have limitations.
Nonetheless, observe that their limitations are opposite to each other: while the first approach is too ``reckless,'' the second approach is too ``conservative.''
This observation motivates us to consider a \emph{combination} of the two approaches which could mitigate each other's limitation. 
\begin{remark}
A similar interpretation of Nesterov's acceleration as a combination of a reckless step and a conservative step also appeared in \citep{allen2014linear,bansal2019potential} 
\end{remark}

Let us implement this idea by alternating between the two approximations \eqref{gd} and \eqref{gd:2-1} of PPM. 
The key modification is that for both approximations, we introduce an additional sequence of points $\{y_t\}$ for cost function approximation; i.e., we use the following approximations for the $t$-th iteration:
\begin{align*}
  f(y_t) +\inp{\nabla f(y_t)}{x-y_t}    \leq f(x) \leq f(y_t) +\inp{\nabla f(y_t)}{x-y_t}+\frac{L}{2}\norm{x-y_t}^2\,.
\end{align*} 
Indeed, this modification is crucial. If we just use approximations at $x_t$, the resulting alternation merely concatenates  \eqref{gd} and \eqref{gd:2-1} during each iteration, and the two limitations we discussed in \S\ref{sec:warmup} will remain in the combined approach. In particular, every other step that corresponds to the lower approximation would be still suffer from overshooting for large stepsizes.

Having introduced a separate sequence $\{y_t\}$ for cost approximations,  we consider the following alternation  where during each iteration, we update $x_{t}$ with  \eqref{gd} and $y_{t}$ with \eqref{gd:2-1}:  
\begin{mdframed}[backgroundcolor=gray!20]
{\bf Approximate PPM with alternating two approaches.} Given $x_0\in \re^d$, let $y_0=x_0$ and run:
\begin{subequations} \label{approx:ppm}
  \begin{align}
    \label{ppm:a} &{\textstyle x_{t+1} \leftarrow  \argmin_{x } \left\{ \low{x}{y_t} + \frac{1}{2\ee_{t+1}} \norm{x-x_t}^2 \right\}}, \\
    \label{ppm:b}
    &{\textstyle y_{t+1}\leftarrow   \argmin_{x} \left\{  \upp{x}{y_t}+ \frac{1}{2\ee_{t+1}} \norm{x-x_{t+1}}^2 \right\} }.
  \end{align} 
\end{subequations} 
\end{mdframed} 
In  Figure~\ref{fig:1}, we empirically verify that \eqref{approx:ppm} indeed gets the best of both worlds: this combined approach successfully approximates PPM even for the regime $\ee_t=\Theta(t)$. 
More remarkably, \eqref{approx:ppm} is exactly equal to one version of AGM (``Form II'' in the introduction).
Turning \eqref{approx:ppm} into the equational form by writing the optimality conditions, and introducing an auxiliary iterate $z_{t+1}:= y_t -\nicefrac{1}{L}\nabla f(y_t)$ (only for simplicity),  we obtain the following ($x_0=y_0=z_0$):
 \begin{mdframed}
 \begin{minipage}{0.50\textwidth} 
 \vspace{-10pt}
 {\bf Equivalent representation of \eqref{approx:ppm}:}
 \begin{subequations} \label{agm}
  \begin{align} 
 & \textstyle y_{t} =\frac{\nicefrac{1}{L}}{ \nicefrac{1}{L}+\ee_t} x_t+\frac{ \ee_t }{  \nicefrac{1}{L}+\ee_t} z_t \,, \label{agm:a}\\
& \textstyle x_{t+1} = x_t -\ee_{t+1} \nabla f(y_{t})\,,\label{agm:b}\\
& \textstyle z_{t+1} = y_{t}-\frac{1}{L}\cdot \nabla f(y_{t})\,. \label{agm:c}
  \end{align} 
\end{subequations}
\end{minipage}
\begin{minipage}{0.48\textwidth} 
\vspace{-10pt}

\begin{figure}[H]
    \centering
    \begin{tikzpicture}[scale=0.34]
\coordinate   (a1) at (-2,0);
\coordinate   (a2) at (-2,-4.5*0.6);
\coordinate   (a3) at (-2,-4.5);
\coordinate [label=above:{ $\boldsymbol{x_t}$}] (xt) at (0,0);
\coordinate [label=above:{ $\boldsymbol{x_{t+1}}$}] (xt1) at (11,0);
\coordinate [label=left:{  $\boldsymbol{z_{t}}$}] (zt) at (2,-4.5);
\coordinate [label=left: {$\boldsymbol{y_{t}}$}] (yt1) at (2*0.6,-4.5*0.6);
\coordinate [label=right: {  $\boldsymbol{z_{t+1}}$}] (zt1) at (6.0,-4.5*0.6);
\draw [->,>=stealth,dotted, line width=1pt,red] (xt) -- node[below] {\tiny $\textcolor{red}{\boldsymbol{-\ee_{t+1} \nabla f(y_{t})}}$} (xt1);
\draw [->,>=stealth, dotted, line width=1pt,red] (yt1) -- node[below] {\tiny $\textcolor{red}{\boldsymbol{-\frac{1}{L}\nabla f(y_{t})}}$} (zt1);

\draw [<->,dotted] (a1) -- node[left] {\small $  {\ee_t}$}  (a2);
\draw [<->,dotted] (a2) -- node[left] {\small $  {\nicefrac{1}{L}}$}  (a3);
\draw [dashed] (zt) --   (zt1);
\draw [dashed] (zt1) --   (xt1);
\draw [dashed] (xt) --   (zt); 
\end{tikzpicture} 
    \caption{Illustration of \eqref{agm}.}
    \label{fig:2}
\end{figure} 

\end{minipage}
 \end{mdframed} 
Hence, we arrive at AGM without relying on any non-trivial derivations in the literature such as estimate sequence~\citep{nesterov2018lectures} or linear coupling~\citep{allen2014linear}.
To summarize, we have demonstrated: 
\begin{center}
    \emph{Nesterov's AGM is an approximate instantiation of the proximal point method!}
\end{center}

\subsection{Understanding mysterious parameters of AGM}
\label{rmk:par}
It is often the case in the literature that the interpolation step \eqref{agm:a} is written as an abstract form $y_t = \tau_t x_t + (1-\tau_t) z_t$ with a weight parameter $\tau_t>0$ to be chosen~\citep{allen2014linear,lessard2016analysis,wilson2016lyapunov,bansal2019potential}.
That said, in the previous works, $\tau_t$ is carefully chosen according to the analysis without conveying  much intuition.  
One important aspect of our PPM view is that it reveals a close relation between the weight parameter $\tau_t$ and the stepsize $\ee_t$.
More specifically, $\tau_t$ is chosen so that the ratio of the distances  $\norm{y_t-x_t}:\norm{y_t-z_t}$ is equal to $\ee_t:\nicefrac{1}{L}$ (see Figure~\ref{fig:2}).

\subsection{Analysis based on PPM perspective}  
\label{subsec:conv}
In order to determine $\ee_t$'s in \eqref{agm},  we revisit the analysis of PPM from \S\ref{sec:warmup}.
In turns out that following \S\ref{sec:app1}, one can derive  the following analogues of \eqref{ineq:1}  and \eqref{ineq:2} using Proposition~\ref{prop:per}  (we defer the derivations to \S\ref{app:agm}):
\begin{align} 
 \tag{\text{$\mathsf{Ineq}^{\mathsf{AGM}}_1$}} &\ee_{t+1}(f(z_{t+1})- f(\xs)) + \frac{1}{2} \norm{\xs -x_{t+1}}^2-
   \frac{1}{2} \norm{\xs -x_{t}}^2\leq \ffo\,,\label{ineq:1:agm}\\
 \tag{\text{$\mathsf{Ineq}^{\mathsf{AGM}}_2$}}&  f(z_{t+1}) -f(z_t) \leq \fft\,,\label{ineq:2:agm}
 \end{align}
 where $\ffo:=(\frac{\ee_{t+1}^2}{2}-\frac{\ee_{t+1}}{2L}) \norm{\nabla f(y_{t})}^2 + L\ee_t\ee_{t+1}\inp{\nabla f(y_{t})}{ z_{t}-y_t}$ and $\fft:= -\frac{1}{2L} \norm{\nabla f(y_{t})}^2- \inp{\nabla f(y_{t})}{z_{t}-y_t}$.
Given the above inequalities, consider the following modified Lyapunov function \eqref{def:lya} which replaces the first $x_t$ with $z_t$: 
\begin{align}
    \textstyle \Phi_t:= (\sum_{i=1}^{t} \ee_i)\cdot (f(z_t)-f(\xs)) + \frac{1}{2}\norm{\xs-x_t}^2\,. \label{lya2}
\end{align}
We note that  \eqref{lya2} is not new; it also appears in prior works~\citep{wilson2016lyapunov,diakonikolas2019approximate,bansal2019potential}, although with different motivations.

Then as before, to prove the validity of the chosen Lyapunov function, it suffices to verify $\ffo + ( \sum_{i=1}^t \ee_i)\cdot \fft \leq 0$, which is equivalent to 
\begin{align}
    \label{cond:1} \textstyle \frac{1}{2L}\left( L\ee_{t+1}^2 - \sum_{i=1}^{t+1} \ee_i\right)   \norm{\nabla f(y_t)}^2 + \left( L\ee_t\ee_{t+1} - \sum_{i=1}^t \ee_i\right)\inp{\nabla f(y_{t})}{z_{t}-y_t} \leq 0
\end{align}
From \eqref{cond:1}, it suffices to choose $\{\ee_t\}$ so that 
$L\ee_t\ee_{t+1}=\sum_{i=1}^t \ee_i$.
Indeed, with such a choice, the coefficient of the inner product term in \eqref{cond:1} becomes zero and the coefficient of the squared norm term becomes $\nicefrac{1}{2L}(L\ee_{t+1}^2 - L\ee_{t+1}\ee_{t+2})\leq 0$ (if $\{\ee_t\}$ is increasing).
Indeed, one can quickly notice that choosing $\ee_{t} =\nicefrac{t}{2L}$ satisfies  the desired relation.
Therefore,  we obtain the well known accelerated convergence rate of $f(z_T)-f(\xs)\leq \frac{2L \norm{x_0-\xs}^2}{T(T+1)}= O(\nicefrac{1}{T^2})$.

 \section{Similar triangle approximations and other variants} 
\label{sec:tri}
In \S\ref{sec:deriv}, we have demonstrated that AGM is nothing but an approximation of PPM.
This view point has not only provided  simple derivations of versions of  AGM, but also offered clear explanations of the stepsizes. In this section, we demonstrate that these interpretations offered by  PPM  actually lead to a great simplification of Nesterov's AGM in the form of the \emph{method of similar triangles}~\citep{nesterov2018lectures,gasnikov2018universal}. 

Our starting point is the observations made in the previous section:  (i) from \S\ref{rmk:par}, we have seen $\norm{y_t-x_t}:\norm{y_t-z_t}=\ee_t:\nicefrac{1}{L}$; (ii) from \S\ref{subsec:conv}, we have seen that we need to choose $\ee_t=\Theta(t)$,  in which case $\ee_{t+1}\approx\ee_t \gg 1$. 
 From these observations, one can readily see that the triangle $\triangle x_tx_{t+1}z_t$ is approximately similar to $\triangle y_tz_{t+1} z_t$.
 Therefore, one can simplify AGM by further exploiting this fact: we  modify the updates so that the two triangles are indeed \emph{similar}.
 There are two different ways one can keep the two triangles similar:
 \begin{enumerate}
 \item We modify the update of $x_{t+1}$ so that the two triangles are similar.
     \item We modify the update of $z_{t+1}$ so that the two triangles are similar.
 \end{enumerate}
 We discuss the above two ways in turn.

\subsection{First similar triangles approximation: momentum form of AGM}
\label{sec:momentum}
We first adopt the first way to keep the two triangles similar. We have the following update.
\begin{mdframed}
 \vspace{-10pt}
 
 \begin{minipage}{0.49\textwidth}
 {\bf First similar triangle approximation:}
 \begin{subequations}\label{sim2}
  \begin{align} 
  & \textstyle y_{t} = \frac{\nicefrac{1}{L}}{ \nicefrac{1}{L}+\ee_t} x_t+\frac{ \ee_t }{ \nicefrac{1}{L}+\ee_t} z_t\,, \label{sim2:a}\\
& \textstyle z_{t+1} = y_t -\frac{1}{L} \nabla f(y_{t})\,,\label{sim2:b}\\
 & \textstyle x_{t+1} =z_{t+1} + L\eta_t (z_{t+1}-z_t)\,. \label{sim2:c}
  \end{align} 
\end{subequations} 
\end{minipage}
\begin{minipage}{0.50\textwidth} 

\begin{figure}[H]
    \centering
    \begin{tikzpicture}[scale=0.35]
\coordinate   (a1) at (-2,0);
\coordinate   (a2) at (-2,-4.5*0.6);
\coordinate   (a3) at (-2,-4.5);
\coordinate [label=above:{ $\boldsymbol{x_t}$}] (xt) at (0,0);
\coordinate [label=above:{ $\boldsymbol{x_{t+1}}$}] (xt1) at (11,0);

\filldraw (2,-4.5) circle (3pt) node[] {};
\coordinate [label=left:{  $\boldsymbol{z_{t}}$}] (zt) at (2,-4.5);

\filldraw (2*0.6,-4.5*0.6) circle (3pt) node[] {};
\coordinate [label=left: {$\boldsymbol{y_{t}}$}] (yt) at (2*0.6,-4.5*0.6);

\filldraw (2+9*0.4,-4.5*0.6) circle (3pt) node[] {};
\coordinate [label=below: {  $\boldsymbol{z_{t+1}}$}] (zt1) at (2+9*0.4,-4.5*0.6);

\filldraw (2+9*0.6,-4.5*0.4) circle (3pt) node[] {};
\coordinate [label=above: {  $\boldsymbol{y_{t+1}}$}] (yt1) at (2+9*0.6,-4.5*0.4);
\draw [->,>=stealth, dotted, line width=0.8pt,red] (xt) -- (xt1);
\draw [->,>=stealth,dotted, line width=0.8pt,red] (yt) --  
node[above] {\tiny $\textcolor{red}{\boldsymbol{-\frac{1}{L} \nabla f(y_{t})}}$} (zt1);

\draw [<->,dotted] (a1) -- node[left] {\small $  {\ee_t}$}  (a2);
\draw [<->,dotted] (a2) -- node[left] {\small $  {\nicefrac{1}{L}}$}  (a3);
\draw [dashed] (zt) --   (xt1);
\draw [dashed] (xt) --   (zt); 
 \draw   pic[draw=black, fill=gray!10, -, angle eccentricity=1.2, angle radius=0.5cm, line width=1pt]
    {angle=xt--xt1--zt};
 \draw   pic[draw=black, fill=gray!10,  -, angle eccentricity=1.2, angle radius=0.5cm, line width=1pt]
    {angle=yt--zt1--zt};
 \draw   pic[draw=black,  fill=gray!10,-, densely dotted, angle eccentricity=1.2, angle radius=0.25cm, line width=1pt]
    {angle=zt--yt--zt1};
 \draw   pic[draw=black, fill=gray!10, -, densely dotted, angle eccentricity=1.2, angle radius=0.25cm, line width=1pt]
    {angle=zt--xt--xt1};
\end{tikzpicture} 
    \caption{The updates of \eqref{sim2}.}
    \label{fig:sim2}
\end{figure}

\end{minipage}
 \end{mdframed} 
 
 In fact, \eqref{sim2} can be equivalently expressed without $\{x_t\}$, as illustrated with dots in Figure~\ref{fig:sim2}.
 More specifically, during the $t$-th iteration, once we compute \eqref{sim2:b}, one can directly update $y_{t+1}$ via $y_{t+1}= z_{t+1} + \frac{L\ee_t}{L\ee_{t+1}+1}(z_{t+1}-z_t)$.
 In other words,
 \begin{align*}   \eqref{sim2} \quad\Longleftrightarrow\quad \begin{cases} 
& \textstyle z_{t+1} = y_{t}-\frac{1}{L}\nabla f(y_{t})\,, \\
& \textstyle y_{t+1} = z_{t+1} +\frac{L\ee_t}{L\ee_{t+1}+1}(z_{t+1}-z_t)\,.
    \end{cases}
\end{align*} 
Hence, \eqref{sim2} is equivalent to the well-known momentum form of AGM (``Form I'' in the introduction).

\paragraph{Recovering popular  stepsize choices.} Notably, our  PPM-based analysis suggests the choice of $\{\ee_t\}$ as per the recursive relation $(L\ee_{t+1}+\frac{1}{2})^2= (L\ee_t+1)^2 +\frac{1}{4}$, which after substitution $L\ee_t+1\leftarrow a_t$ exactly recovers the popular recursive relation $a_{t+1} = \frac{1}{2} (1+\sqrt{1+4a_t^2})$ in \citep{nesterov1983method,beck2009fast}. 
 The analysis is similar to the one given in \S\ref{subsec:conv}. Below we provide the details.

 Following \S\ref{subsec:conv}, we again derive the following  counterparts of \eqref{ineq:1} and \eqref{ineq:2} with straightforward arguments (see  \S\ref{app:rmk} for details):
\begin{align}
&\tag{\text{$\mathsf{Ineq}^{\mathsf{SIM}}_1$}} \eep_{t+1}[f(z_{t+1})-f(\xs)] +\frac{1}{2}\norm{\xs-x_{t+1}}^2  -\frac{1}{2}\norm{\xs-x_t}^2  \leq \hho\,, \label{ineq:1:mom}\\
&\tag{\text{$\mathsf{Ineq}^{\mathsf{SIM}}_2$}} f(z_{t+1}) -f(z_t) \leq \hht \label{ineq:2:mom}\,,
\end{align}
where using the notation $\eep_{t+1} :=\ee_t+\frac{1}{L}$, the right hand side of the above inequalites are defined as $\hho:=\frac{1}{2}\left(  -(L\ee_t+1)^2   +L\eep_{t+1}\right)\cdot \norm{z_{t+1}-y_t}^2 +  \eep_{t+1}\cdot \inp{\nabla f(y_{t})}{z_{t+1}-x_{t+1}}$ and $\hht:=\frac{L}{2} \norm{z_{t+1}-y_t}^2+ \inp{\nabla f(y_{t})}{z_{t+1}-z_t}$.

Having established counterparts of \eqref{ineq:1} and \eqref{ineq:2}, following \S\ref{subsec:conv}, we choose
\begin{align}
    \textstyle \Phi_t:= (\sum_{i=1}^{t} \eep_i)\cdot (f(z_t)-f(\xs)) + \frac{1}{2}\norm{\xs-x_t}^2\,. \label{lya2:2}
\end{align}
To prove the validity of the chosen Lyapunov function, it suffices to verify 
\begin{align} \label{suff:2}
\textstyle    \hho + ( \sum_{i=1}^t \eep_i)\cdot \hht \leq 0
\end{align}
which is equivalent to showing (because $z_{t+1}-x_{t+1} = -L\ee_t (z_{t+1}-z_t) $):
\begin{align}
    \label{cond:1:2}
    \begin{split}
        \textstyle \frac{1}{2}\left(  -(L\ee_t+1)^2     +\sum_{i=1}^{t+1} L\eep_i \right)\cdot \norm{z_{t+1}-y_t}^2  + \textstyle \left( L\ee_t\eep_{t+1} - \sum_{i=1}^t \eep_i\right)\inp{\nabla f(y_{t})}{z_{t+
    1}-z_t} 
    \end{split} \leq 0\,.
\end{align}
From \eqref{cond:1:2}, it suffices to choose $\{\ee_t\}$ so that $L\ee_t\eep_{t+1}=\sum_{i=1}^t \eep_i$.
Indeed, with such a choice, the coefficient of the inner product term in \eqref{cond:1:2} becomes zero and the coefficient of the squared norm term becomes
\begin{align*}
   \textstyle \frac{1}{2}\left(  -(L\ee_t+1)^2     +\sum_{i=1}^{t+1} L\eep_i \right)&= \frac{1}{2}\left(  -(L\ee_t+1)^2     +L\eep_{t+1}+ L\eep_{t+1}\cdot L\ee_t\right) \\
    &\textstyle= \frac{1}{2}\left(  -(L\ee_t+1)^2     +L\eep_{t+1} (L\ee_t+1) \right) =0
\end{align*}
since $L\eep_{t+1} = L\ee_t+1$. 
Indeed, one can actually simplify the relation $L\ee_t\eep_{t+1}=\sum_{i=1}^t \eep_i$:
\begin{align*}
   \textstyle L\ee_{t+1} \cdot (L\ee_{t+1}+1)= L\ee_{t+1}\cdot L\eep_{t+2}=\sum_{i=1}^{t+1} L\eep_i  = L\eep_{t+1} + L\ee_{t}\cdot L\eep_{t+1} = (L\ee_{t}+1)^2\,.
\end{align*}
After rearranging, we obtain the recursive relation:
$(L\ee_{t+1}+\frac{1}{2})^2= (L\ee_t+1)^2 +\frac{1}{4}$, which after the substitution $L\ee_t+1= a_t$ exactly recovers the popular recursive relation $a_{t+1} = \frac{1+\sqrt{1+4a_t^2}}{2}$ in \citep{nesterov1983method,beck2009fast}.

 \subsection{Second similar triangles approximation: acceleration for composite costs}
 \label{sec:composite}

 We now adopt the second way to keep the two triangles similar. We have the following update.
 
\begin{mdframed}

 \vspace{-10pt}
 \begin{minipage}{0.49\textwidth} 
  {\bf Second similar triangle approximation:}
 \begin{subequations}\label{sim}
  \begin{align} 
  & \textstyle y_{t} = \frac{\nicefrac{1}{L}}{ \nicefrac{1}{L}+\ee_t} x_t+\frac{ \ee_t }{ \nicefrac{1}{L}+\ee_t} z_t\,, \label{sim:a}\\
& \textstyle x_{t+1} = x_t -\ee_{t+1} \nabla f(y_{t})\,,\label{sim:b}\\
 & \textstyle z_{t+1} = \frac{\nicefrac{1}{L}}{ \nicefrac{1}{L}+\ee_t} x_{t+1}+\frac{ \ee_t }{ \nicefrac{1}{L}+\ee_t } z_t\,. \label{sim:c}
  \end{align} 
\end{subequations} 
\end{minipage}
\begin{minipage}{0.50\textwidth} 
\vspace{-5pt}
\begin{figure}[H]
    \centering
    \begin{tikzpicture}[scale=0.34]
\coordinate   (a1) at (-2,0);
\coordinate   (a2) at (-2,-4.5*0.6);
\coordinate   (a3) at (-2,-4.5);
\coordinate [label=above:{ $\boldsymbol{x_t}$}] (xt) at (0,0);
\coordinate [label=above:{ $\boldsymbol{x_{t+1}}$}] (xt1) at (11,0);
\coordinate [label=left:{  $\boldsymbol{z_{t}}$}] (zt) at (2,-4.5);
\coordinate [label=left: {$\boldsymbol{y_{t}}$}] (yt) at (2*0.6,-4.5*0.6);
\coordinate [label=right: {  $\boldsymbol{z_{t+1}}$}] (zt1) at (2+9*0.4,-4.5*0.6);
\draw [->,>=stealth, dotted, line width=0.8pt,red] (xt) -- node[below] {\tiny $\textcolor{red}{\boldsymbol{-\ee_{t+1} \nabla f(y_{t})}}$} (xt1);
\draw [->,>=stealth,dotted, line width=0.8pt,red] (yt) --  (zt1);

\draw [<->,dotted] (a1) -- node[left] {\small $  {\ee_t}$}  (a2);
\draw [<->,dotted] (a2) -- node[left] {\small $  {\nicefrac{1}{L}}$}  (a3);
\draw [dashed] (zt) --   (xt1);
\draw [dashed] (xt) --   (zt); 
 \draw   pic[draw=black, fill=gray!10, -, angle eccentricity=1.2, angle radius=0.5cm, line width=1pt]
    {angle=xt--xt1--zt};
 \draw   pic[draw=black, fill=gray!10,  -, angle eccentricity=1.2, angle radius=0.5cm, line width=1pt]
    {angle=yt--zt1--zt};
 \draw   pic[draw=black,  fill=gray!10,-, densely dotted, angle eccentricity=1.2, angle radius=0.25cm, line width=1pt]
    {angle=zt--yt--zt1};
 \draw   pic[draw=black, fill=gray!10, -, densely dotted, angle eccentricity=1.2, angle radius=0.25cm, line width=1pt]
    {angle=zt--xt--xt1};
\end{tikzpicture} 
    \caption{Illustration of \eqref{sim}.}
    \label{fig:sim}
\end{figure}
\end{minipage}
 \end{mdframed}
 \noindent This is ``Form III'' in the introduction.
 Below, we provide a PPM-based analysis  for a more general setting.

    One advantage of   \eqref{sim} is that it  admits a simple extension to the practical  setting of constrained optimization on composite costs (see e.g.  \citep[\S6.1.3]{nesterov2018lectures} for applications).
  More specifically, for a closed convex set $Q\subseteq \re^d$ and a closed\footnote{This means that the epigraph of the function is closed. See \citep[Definition 3.1.2]{nesterov2018lectures}. } convex function $\Psi:Q\to \re$, consider \begin{align*}
    \textstyle \min_{x\in Q} \fsi(x):= f(x)+\Psi(x) \,,
\end{align*}
where $f:Q\to \re$ is a differentiable convex function which is $L$-smooth with respect to a norm $\norm{\cdot}$ that is not necessarily  the $\ell_2$ norm (i.e., we regard the norm in Definition~\ref{def:ell2} to be our chosen norm).  For the general norm case, we use the Bregman divergence.
\begin{definition}
Given a $1$-strongly convex (w.r.t the chosen norm $\norm{\cdot}$) function $h:Q\to \re \cup \{\infty\}$  that is differentiable on the interior of $Q$, $\breg{h}{u}{v}:=h(u)-h(v)-\inp{\nabla h(v)}{u-v}$ for all $u,v\in Q$.
\end{definition} 
\noindent Under the above setting and assumption, \eqref{sim} admits a simple generalization:
\begin{mdframed}
 {\bf Generalization of \eqref{sim} to composite costs:}
 \begin{subequations}\label{simg}
  \begin{align} 
 &\textstyle y_{t} =\frac{\nicefrac{1}{L}}{\nicefrac{1}{L}+\ee_t } x_t+\frac{ \ee_t }{  \nicefrac{1}{L}+\ee_t } z_t\,, \label{simg:a}\\
&\textstyle x_{t+1} \leftarrow \argmin_{x\in Q} \left\{ \low{x}{y_t}+ \frac{1}{\ee_{t+1}} \breg{h}{x}{x_t} +\Psi(x) \right\}\,,\label{simg:b}\\
& \textstyle z_{t+1} =\frac{\nicefrac{1}{L}}{\nicefrac{1}{L}+\ee_t} x_{t+1}+\frac{ \ee_t }{\nicefrac{1}{L}+ \ee_t } z_t \,. \label{simg:c}
  \end{align} 
\end{subequations}  
\end{mdframed} 
 
  Now we provide a simple PPM-based analysis of \eqref{simg}:
\begin{proof}[{\bf PPM-based analysis of \eqref{simg}}] 
To obtain counterparts of \eqref{ineq:1} and \eqref{ineq:2}, we now use a generalization of Proposition \ref{prop:per} to the Bregman divergence \citep[Lemma 3.1]{teboulle2018simplified}. 
With such a generalization,  we obtain the following inequality for  $\pss(x):=\ee_{t+1}[f(y_t)+\inp{\nabla f(y_t)}{x-y_t} +\Psi(x)]$:
\begin{align}
    \pss(x_{t+1}) -\pss(\xs) +  \breg{h}{\xs}{x_{t+1}}+\breg{h}{x_{t+1}}{x_{t}}-
  \breg{h}{\xs}{x_{t}} \leq 0\,, \label{gen:1}
\end{align}
where $\xs\in \argmin_{x\in Q}\fsi(x)$.
Now using \eqref{gen:1}, one can derive from first principles the following inequalities   (we defer the derivations to \S\ref{app:sim}):
\begin{align}
   \tag{\text{$\mathsf{Ineq}^{\mathsf{SIM'}}_1$}}&\ee_{t+1}(\fsi(z_{t+1})- \fsi(\xs) ) + \breg{h}{\xs}{x_{t+1}}-
  \breg{h}{\xs}{x_{t}} \leq \ggo\,, \label{ineq:1:sim}\\
    \tag{\text{$\mathsf{Ineq}^{\mathsf{SIM'}}_2$}}  &\fsi(z_{t+1}) -\fsi(z_t) \leq \ggt\,. \label{ineq:2:sim}
    \end{align}
where $\ggo:=   -\frac{1}{2}\norm{ x_{t+1}-x_t}^2  +\ee_{t+1}[\frac{L}{2}\norm{z_{t+1}-y_{t}}^2 +\inp{\nabla f(y_{t})}{z_{t+1}-x_{t+1}}   +\Psi(z_{t+1}) -\Psi(x_{t+1})]$ and $\ggt:= \frac{L}{2} \norm{z_{t+1}-y_t}^2 +  \inp{\nabla f(y_{t})}{z_{t+1}-z_t} +\Psi(z_{t+1})-\Psi(z_{t}) \nonumber$.
    Similar to \S\ref{subsec:conv}, yet replacing the norm squared term with the Bregman divergence, we choose
    \begin{align*}
\textstyle\Phi_t:= (\sum_{i=1}^{t} \ee_i)\cdot (\fsi(z_t)-\fsi(\xs)) +\breg{h}{\xs}{x_t}.
    \end{align*}
    Then, it suffices to show  $\ggo + ( \sum_{i=1}^t \ee_i)\cdot \ggt \leq 0$.
    Using   the facts (i) $z_{t+1}-x_{t+1} = L\ee_t (z_t-z_{t+1})$ and (ii) $\norm{x_{t+1}-x_t}= (L\ee_t+1)\norm{z_{t+1}-y_t}$  (both are immediate consequences of the similar triangles) and rearranging, one can easily check that  $\ggo + ( \sum_{i=1}^t \ee_i)\cdot \ggt$ is equal to  
    \begin{align}
\textstyle     \frac{1}{2}\left(-(L\ee_t+1)^2 +L \ee_{t+1}+  L\sum_{i=1}^{t}\ee_{i} \right)\norm{z_{t+1}-y_{t}}^2 \label{a1}\\
\textstyle +\left(L\ee_t \ee_{t+1} - \sum_{i=1}^t\ee_i\right)\inp{\nabla f(y_{t})}{z_{t}-z_{t+1}} \label{a2}\\
\textstyle  +\ee_{t+1}[ \Psi(z_{t+1}) -\Psi(x_{t+1}) ] +   \left( \sum_{i=1}^t \ee_i\right)\cdot [\Psi(z_{t+1})-\Psi(z_{t})].
   \label{a3}
\end{align}
    Now choosing $\ee_t= \nicefrac{t}{2L}$ analogously to \S\ref{subsec:conv}, one can easily verify  $\eqref{a1}+\eqref{a2}+\eqref{a3} \leq 0$.
    Indeed, for \eqref{a1}, since  $L  \ee_t \ee_{t+1}=\sum_{i=1}^t\ee_i$, the coefficient becomes $\nicefrac{1}{2}(L\ee_t+1)(L\ee_{t+1}-L\ee_t-1)$ which is a negative number since $L\ee_{t+1}-L\ee_t-1=-\nicefrac{1}{2}$; 
    for \eqref{a2}, the coefficient becomes zero due to the relation $L  \ee_t \ee_{t+1}=\sum_{i=1}^t\ee_i$; 
    lastly, for \eqref{a3}, we have
    \begin{align} 
\eqref{a3}=    \ee_{t+1}\left[ (1+L\ee_t) \Psi(z_{t+1}) -\Psi(x_{t+1}) -L\ee_{t}\Psi(z_{t})\right] \leq 0\,,
    \end{align}
    where the equality is due to the relation $L  \ee_t \ee_{t+1} =\sum_{i=1}^t\ee_i$, and the inequality is due to the update \eqref{simg:c} (which can be equivalently written as $(1+L\ee_t) z_{t+1} = x_{t+1}+ L\ee_t z_t$) and the convexity of $\Psi$.
    Hence, we obtain the accelerated rate of  $\fsi(z_T)-\fsi(\xs)\leq \frac{4L\breg{h}{\xs}{x_0}}{T(T+1)} = O(\nicefrac{1}{T^2})$.
\end{proof}

\section{Extension to strongly convex costs} 
\label{sec:str}

In this section, we  extend our PPM framework to the case of strongly convex costs. 
As we shall see, our framework gives rise to a simple derivation of the most general version of AGM  called  “\emph{General Scheme for Optimal Method}”~\citep[(2.2.7)]{nesterov2018lectures}.  
We first make  the approximate PPM \eqref{approx:ppm} more flexible by considering two separate stepsizes.
\begin{mdframed} 
{\bf Approximate PPM with two separate stepsizes $\{\ee_t\}$ and $\{\eep_t\}$.} Given $x_0=y_0\in \re^d$, 
\begin{subequations}\label{approx:ppm2}
 \begin{align}
    \label{ppm2:a} &{\textstyle x_{t+1} \leftarrow  \argmin_{x } \left\{ \low{x}{y_t} + \frac{1}{2\ee_{t+1}} \norm{x-x_t}^2 \right\}}, \\
    \label{ppm2:b}
    &{\textstyle y_{t+1}\leftarrow   \argmin_{x} \left\{  \upp{x}{y_t}+ \frac{1}{2\eep_{t+1}} \norm{x-x_{t+1}}^2 \right\} }.
  \end{align} 
\end{subequations} 
\end{mdframed} 
Now let us  apply our PPM view  to the strongly convex cost case. 
\begin{definition}[Strong convexity] \label{def:str}
For $\mu>0$, we say a differentiable function $f:\re^d\to \re$ is  $\mu$-strongly convex if   $f(x) \geq f(y) + \inp{\nabla f(y)}{x-y} +\frac{\mu}{2}\norm{x-y}^2 $ for any $x,y\in\re^d$.
\end{definition} 
\noindent Since $f$ is additionally assumed to be strongly convex, one can now strengthen the lower approximation $\low{x}{y_t}$ in \eqref{ppm2:a} to $\low{x}{y_t} + \frac{\mu}{2}\norm{x-y_t}^2$.
In other words, we obtain
\begin{mdframed}
{\bf Approximate PPM for strongly-convex costs.} Given $x_0=y_0\in \re^d$,
\begin{subequations}\label{approx:ppm:strong}
 \begin{align}
    \label{ppm:strong:a} & x_{t+1} \leftarrow \argmin_{x\in \re^d} \bigg\{ \low{x}{y_t} +\underbrace{\frac{\mu}{2}\norm{x-y_t}^2}_{ \substack{\text{additional term due to}\\ \text{strong convexity}}} + \frac{1}{2\ee_{t+1}} \norm{x-x_t}^2 \bigg\}, \\
    \label{ppm:strong:b}
    &{\textstyle y_{t+1}\leftarrow   \argmin_{x} \left\{  \upp{x}{y_t}+ \frac{1}{2\eep_{t+1}} \norm{x-x_{t+1}}^2 \right\} }.
  \end{align} 
\end{subequations} 
\end{mdframed} 
Writing the optimality condition of \eqref{approx:ppm:strong}, it is straightforward to check that the  approximate PPM \eqref{approx:ppm2} is equivalent to the following updates ($x_0=y_0=z_0$): 
 \begin{mdframed}
  
  \vspace{-10pt}
  
 \begin{minipage}{0.5\textwidth}
 {\bf Equivalent representation of \eqref{approx:ppm:strong}:}
 \begin{subequations}\label{agm2}
  \begin{align} 
  & \textstyle y_{t} =\frac{\nicefrac{1}{L}}{  \nicefrac{1}{L}+\eep_t} x_t +\frac{ \eep_t }{   \nicefrac{1}{L}+\eep_t} z_t\,, \label{agm2:a}\\ 
 \begin{split}
    &\textstyle x_{t+1} = \frac{\nicefrac{1}{\mu}}{ \nicefrac{1}{\mu}+\ee_{t+1}}x_t +\frac{\ee_{t+1}}{ \nicefrac{1}{\mu} +\ee_{t+1}}y_t\\
    &\textstyle \quad \quad - \frac{\nicefrac{1}{\mu}\cdot \ee_{t+1}}{ \nicefrac{1}{\mu}+\ee_{t+1}}\nabla f(y_{t})
\end{split}\,,\label{agm2:b}\\
& \textstyle z_{t+1} = y_{t}-\frac{1}{L}\nabla f(y_{t})\,. \label{agm2:c}
  \end{align} 
\end{subequations} 

\end{minipage} 
\begin{minipage}{0.5\textwidth}

\begin{figure}[H]
    \centering
    \begin{tikzpicture}[scale=0.33]
\coordinate   (a1) at (-3,0);
\coordinate   (a2) at (-3,-7*0.7);
\coordinate   (a3) at (-3,-7);

\coordinate   (b1) at (-2.5,0);
\coordinate   (b2) at (-2.5,-7*0.3);
\coordinate   (b3) at (-2.5,-7*0.7);
\coordinate [label=above:{ $\boldsymbol{x_t}$}] (xt) at (0,0);
\coordinate [label=above:{$\boldsymbol{x_{t+1}}$}] (xt1) at (11,0);
\coordinate  [label=left:{ $\boldsymbol{w_t}$}]   (wt) at (1*0.3,-7*0.3);
\coordinate [label=left:{  $\boldsymbol{z_{t}}$}] (zt) at (1,-7);
\coordinate [label=left: {$\boldsymbol{y_{t}}$}] (yt1) at (1*0.7,-7*0.7);
\coordinate [label=right: {$\boldsymbol{z_{t+1}}$}] (zt1) at (1+10*0.3/0.7,-7+7*0.3/0.7);
\draw [->,>=stealth,dotted, line width=1pt,red] (wt) -- node[above] {\tiny $\textcolor{red}{\boldsymbol{-\frac{\ee_{t+1}\cdot \nicefrac{1}{\mu}}{ \ee_{t+1}+\nicefrac{1}{\mu}} \nabla f(y_{t})}}$} (xt1);
\draw [->,>=stealth, dotted, line width=1pt,red] (yt1) -- node[above] {\tiny $\textcolor{red}{\boldsymbol{-\frac{1}{L}\nabla f(y_{t})}}$} (zt1);

\draw [<->,dotted] (a1) -- node[left] {\small $  {\eep_t}$}  (a2);
\draw [<->,dotted] (a2) -- node[left] {\small $  {\nicefrac{1}{L}}$}  (a3);
\draw [dashed] (zt) --   (zt1);
\draw [dashed] (zt1) --   (xt1);
\draw [dashed] (xt) --   (zt); 

\draw [<->,dotted] (b1) -- node[right] {\scriptsize $  {\ee_{t+1}}$}  (b2);
\draw [<->,dotted] (b2) -- node[right] {\small $  {\nicefrac{1}{\mu}}$}  (b3); 
\end{tikzpicture} 
    \caption{Illustration of \eqref{agm2}.}
    \label{fig:2-2}
\end{figure}  

\end{minipage}
 \end{mdframed} 
 Note that \eqref{agm2} is the most general version of AGM due to Nesterov called  “\emph{General Scheme for Optimal Method}”~\citep[(2.2.7)]{nesterov2018lectures} (``Form IV'' in the introduction).
Again, our derivation provides new insights into the choices of the AGM stepsizes by expressing them in terms of the PPM stepsizes $\ee_t$'s and $\eep_t$'s.

\subsection{Relation to well known momentum version}
\label{sec:str:momen}
Perhaps, the most well known version of AGM for strongly convex costs is the  momentum version due to Nesterov (see, e.g., \citep[(2.2.22)]{nesterov2018lectures})
\begin{align} \label{agm:str:nesterov}
\begin{split}
& \textstyle z_{t+1} = y_{t}-\frac{1}{L}\nabla f(y_{t})\,, \\
& \textstyle y_{t+1} = z_{t+1} +\frac{\sqrt{\kappa}-1}{\sqrt{\kappa}+1}(z_{t+1}-z_t)\,.    
\end{split}
\end{align}
One might wonder whether one can better understand the stepsizes in \eqref{agm:str:nesterov} from \eqref{agm2}. 

Let us first recall the well known convergence rate of PPM for strongly convex costs due to \citep[(1.14)]{rockafellar1976monotone}:
\begin{align}\label{conv:prox:str}
    \textstyle f(x_{T})- f(\xs) \leq \OO{\prod_{t=1}^{T} (1+\mu\ee_t)^{-1}}\quad \text{for any $T\geq 1$.} 
\end{align}
From \eqref{conv:prox:str}, one can see that in order to achieve the accelerated convergence rate $O(\exp(\nicefrac{-T}{\sqrt{\kappa}}))$ where $\kappa$ is the condition number $\nicefrac{L}{\mu}$,  the stepsizes $\ee_t$ must be chosen so that $\ee_t \approx \mu^{-1}(\sqrt{\kappa})^{-1}$.
In fact, the well known version \eqref{agm:str:nesterov} corresponds to choosing the following stepsizes for \eqref{agm2}:
\begin{align}\label{stepchoice}
    \text{$\ee_t \equiv \ee:= \mu^{-1}  (\sqrt{\kappa}-1)^{-1}$\quad and \quad$\eep_t \equiv \eep:= \mu^{-1}  (\sqrt{\kappa})^{-1}$.} 
\end{align}  
To see this, note that with such choice of $\ee$ and $\eep$, \eqref{agm2} becomes:
 \begin{mdframed}
 \vspace{-10pt}
  
  \begin{minipage}{0.5\textwidth} 
  
  {\bf \eqref{agm2} with stepsize chosen as \eqref{stepchoice}:} 
\begin{subequations}\label{agm2p}
  \begin{align} 
  & \textstyle y_{t} =\frac{1}{ 1+\sqrt{\kappa}} x_t +\frac{\kappa}{ 1+\sqrt{\kappa}} z_t\,, \label{agm2p:a}\\
& \textstyle x_{t+1} = \frac{\sqrt{\kappa}-1}{ \sqrt{\kappa}}x_t +\frac{1}{ \sqrt{\kappa}}y_t- \frac{\sqrt{\kappa}}{L}\nabla f(y_{t})\,,\label{agm2p:b}\\
& \textstyle z_{t+1} = y_{t}-\frac{1}{L}\nabla f(y_{t})\,. \label{agm2p:c}
  \end{align} 
\end{subequations} 
 
\end{minipage}
\begin{minipage}{0.5\textwidth} 
\begin{figure}[H]
    \centering
    \begin{tikzpicture}[scale=0.33]
\coordinate   (a1) at (-3.7,0);
\coordinate   (a2) at (-3.7,-7*0.7);
\coordinate   (a3) at (-3.7,-7);

\coordinate   (b1) at (-3.2,0);
\coordinate   (b2) at (-3.2,-7*0.3);
\coordinate   (b3) at (-3.2,-7*0.7);
\coordinate [label=above:{ $\boldsymbol{x_t}$}] (xt) at (0,0);
\coordinate [label=above:{$\boldsymbol{x_{t+1}}$}] (xt1) at (11,0);
\coordinate  [label=left:{ $\boldsymbol{w_t}$}]   (wt) at (1*0.3,-7*0.3);
\coordinate [label=left:{  $\boldsymbol{z_{t}}$}] (zt) at (1,-7);
\coordinate [label=left: {$\boldsymbol{y_{t}}$}] (yt) at (1*0.7,-7*0.7);
\coordinate [label=right: {$\boldsymbol{z_{t+1}}$}] (zt1) at (1+10*0.3/0.7,-7+7*0.3/0.7);
\draw [->,>=stealth,dotted, line width=1pt,red] (wt) -- node[above] {\tiny $\textcolor{red}{\boldsymbol{-\frac{\sqrt{\kappa}}{L} \nabla f(y_{t})}}$} (xt1);

\draw [->,>=stealth, dotted, line width=1pt,red] (yt) -- node[above] {\tiny $\textcolor{red}{\boldsymbol{-\frac{1}{L}\nabla f(y_{t})}}$} (zt1);

\draw [<->,dotted] (a1) -- node[left] {\small $  {\sqrt{\kappa}}$}  (a2);
\draw [<->,dotted] (a2) -- node[left] {\small $  {1}$}  (a3);
\draw [dashed] (zt) --   (zt1);
\draw [dashed] (zt1) --   (xt1);
\draw [dashed] (xt) --   (zt); 

\draw [<->,dotted] (b1) -- node[right] {\small $  {1}$}  (b2);
\draw [<->,dotted] (b2) -- node[right] {\scriptsize $  {\sqrt{\kappa}-1}$}  (b3); 
 \draw   pic[draw=black, fill=gray!10, -, angle eccentricity=1.2, angle radius=0.5cm, line width=1pt]
    {angle=wt--xt1--zt};
 \draw   pic[draw=black, fill=gray!10,  -, angle eccentricity=1.2, angle radius=0.5cm, line width=1pt]
    {angle=yt--zt1--zt};
     \draw   pic[draw=black,  fill=gray!10,-, densely dotted, angle eccentricity=1.2, angle radius=0.2cm, line width=1pt]
    {angle=zt--yt--zt1};
 \draw   pic[draw=black, fill=gray!10, -, densely dotted, angle eccentricity=1.2, angle radius=0.2cm, line width=1pt]
    {angle=zt--wt--xt1};
\end{tikzpicture} 
    \caption{Illustration of \eqref{agm2p}.}
    \label{fig:2-2p}
\end{figure}   
\end{minipage}
 \end{mdframed}
 As shown in Figure~\ref{fig:2-2p},  $\triangle w_tx_{t+1}z_t$ is  similar to $\triangle y_tz_{t+1} z_t$, so one can write the updates \eqref{agm2p}  without $\{x_t\}$ and $\{w_t\}$, which precisely recovers \eqref{agm:str:nesterov}:
\begin{align*}   \eqref{agm2p} \quad\Longleftrightarrow\quad \eqref{agm:str:nesterov}=\begin{cases} 
& \textstyle z_{t+1} = y_{t}-\frac{1}{L}\nabla f(y_{t})\,, \\
& \textstyle y_{t+1} = z_{t+1} +\frac{\sqrt{\kappa}-1}{\sqrt{\kappa}+1}(z_{t+1}-z_t)\,.
    \end{cases}
\end{align*}

\section{Related work} \label{sec:related}

Our approach is inspired by that of \citet{defazio2019curved} that  establishes an inspiring connection  between AGM and PPM. 
The main observation in that paper is that  for strongly convex costs, one can derive  a version of AGM from the primal-dual form of PPM with a tweak of geometry.
Compared with \citep{defazio2019curved}, our approach strengthens the connection between AGM and PPM by considering more versions of AGM and their analyses.
Another advantage of our approach is that it does not require duality.

We now summarize previous works on developing alternative approaches to Nesterov's acceleration.
Most works have studied  the continuous limit dynamics of Nesterov's AGM~\citep{su2014differential,krichene2015accelerated,wibisono2016variational}.
These continuous dynamics approaches have brought about new intuitions about Nesterov's acceleration, and follow-up works have developed analytical techniques for such dynamics~\citep{wilson2016lyapunov,diakonikolas2019approximate}. 
Another notable contribution is made based on the  linear coupling framework~\citep{allen2014linear}.
The main observation is that the two most popular first-order methods, namely gradient descent and mirror descent, have complementary performances, and hence, one can come up with a faster method by linearly coupling  the two methods.   
Lastly, Nesterov's acceleration has  been explained from the perspective of  computing the equilibrium in a primal-dual game~\citep{wang2018acceleration,cohen2021relative}.

  PPM has been used to design or interpret other  optimization methods~\citep{drusvyatskiy2017proximal}.
 To list few instances, PPM has given rise to fast methods  for weakly convex problems~\citep{davis2019proximally}, the prox-linear methods for composite optimizations~\citep{burke1995gauss,nesterov2007modified,lewis2016proximal}, accelerated methods for stochastic optimizations~\citep{lin2015universal}, and methods for saddle-point problems~\citep{mokhtari2019unified}.

\section{Conclusion}	 
	 
	 This work provides a way to understand Nesterov's acceleration based on the proximal point method.
The framework presented in this paper motivates a simplification of AGM using similar triangles and readily extends to the strongly convex case and recovers the most general accelerated method due to Nesterov.

We believe that the simple derivations presented in this paper clarify and deepen our understanding of Nesterov's acceleration. 
Our framework is therefore not only of pedagogical value but also helpful for research. 
For future directions, it would be interesting to connect our PPM view to accelerated stochastic methods
~\citep{lin2015universal,lan2018optimal} and other accelerated methods, including geometric descent~\citep{bubeck2015geometric}.
Furthermore, we hope the connections presented in this work will help advance the development of accelerated methods in settings much wider than convex optimization (see e.g.,~\citep{bacak2014convex}).

\section*{Acknowledgement} 
We thank {\bf Alp Yurtsever} and   {\bf Jingzhao Zhang}  for detailed comments and stimulating discussions, {\bf Aaron Defazio} for clarifications that help the author develop \S\ref{sec:composite}, and  {\bf Heinz Bauschke} for constructive suggestions on the presentation of \S\ref{sec:deriv} and \S\ref{sec:str:momen}. 
Kwangjun Ahn  and Suvrit Sra acknowledge support from the NSF Grant (CAREER: 1846088).
Kwangjun Ahn also acknowledge support from  Kwanjeong Educational Foundation.

	\bibliographystyle{plainnat}
	\bibliography{ref}

\begin{thebibliography}{36}
\providecommand{\natexlab}[1]{#1}
\providecommand{\url}[1]{\texttt{#1}}
\expandafter\ifx\csname urlstyle\endcsname\relax
  \providecommand{\doi}[1]{doi: #1}\else
  \providecommand{\doi}{doi: \begingroup \urlstyle{rm}\Url}\fi

\bibitem[Allen-Zhu and Orecchia(2017)]{allen2014linear}
Zeyuan Allen-Zhu and Lorenzo Orecchia.
\newblock Linear coupling: An ultimate unification of gradient and mirror
  descent.
\newblock In \emph{ITCS 2017}. Schloss Dagstuhl-Leibniz-Zentrum fuer
  Informatik, 2017.

\bibitem[Auslender and Teboulle(2006)]{auslender2006interior}
Alfred Auslender and Marc Teboulle.
\newblock Interior gradient and proximal methods for convex and conic
  optimization.
\newblock \emph{SIAM Journal on Optimization}, 16\penalty0 (3):\penalty0
  697--725, 2006.

\bibitem[Bac{\'a}k(2014)]{bacak2014convex}
Miroslav Bac{\'a}k.
\newblock \emph{Convex analysis and optimization in Hadamard spaces},
  volume~22.
\newblock Walter de Gruyter GmbH \& Co KG, 2014.

\bibitem[Bansal and Gupta(2019)]{bansal2019potential}
Nikhil Bansal and Anupam Gupta.
\newblock Potential-function proofs for gradient methods.
\newblock \emph{Theory of Computing}, 15\penalty0 (4):\penalty0 1--32, 2019.

\bibitem[Bauschke and Combettes(2011)]{bauschke2011convex}
Heinz~H Bauschke and Patrick~L Combettes.
\newblock \emph{Convex analysis and monotone operator theory in {Hilbert}
  spaces}, volume 408.
\newblock Springer, 2011.

\bibitem[Beck and Teboulle(2009)]{beck2009fast}
Amir Beck and Marc Teboulle.
\newblock A fast iterative shrinkage-thresholding algorithm for linear inverse
  problems.
\newblock \emph{SIAM journal on imaging sciences}, 2\penalty0 (1):\penalty0
  183--202, 2009.

\bibitem[Becker et~al.(2011)Becker, Bobin, and Cand{\`e}s]{becker2011nesta}
Stephen Becker, J{\'e}r{\^o}me Bobin, and Emmanuel~J Cand{\`e}s.
\newblock {NESTA}: A fast and accurate first-order method for sparse recovery.
\newblock \emph{SIAM Journal on Imaging Sciences}, 4\penalty0 (1):\penalty0
  1--39, 2011.

\bibitem[Bubeck et~al.(2015)Bubeck, Lee, and Singh]{bubeck2015geometric}
S{\'e}bastien Bubeck, Yin~Tat Lee, and Mohit Singh.
\newblock A geometric alternative to {Nesterov's} accelerated gradient descent.
\newblock \emph{arXiv preprint: 1506.08187}, 2015.

\bibitem[Burke and Ferris(1995)]{burke1995gauss}
James~V Burke and Michael~C Ferris.
\newblock A {Gauss}-{Newton} method for convex composite optimization.
\newblock \emph{Mathematical Programming}, 71\penalty0 (2):\penalty0 179--194,
  1995.

\bibitem[Cohen et~al.(2021)Cohen, Sidford, and Tian]{cohen2021relative}
Michael~B Cohen, Aaron Sidford, and Kevin Tian.
\newblock Relative lipschitzness in extragradient methods and a direct recipe
  for acceleration.
\newblock In \emph{12th Innovations in Theoretical Computer Science Conference
  (ITCS 2021)}. Schloss Dagstuhl-Leibniz-Zentrum f{\"u}r Informatik, 2021.

\bibitem[Davis and Grimmer(2019)]{davis2019proximally}
Damek Davis and Benjamin Grimmer.
\newblock Proximally guided stochastic subgradient method for nonsmooth,
  nonconvex problems.
\newblock \emph{SIAM Journal on Optimization}, 29\penalty0 (3):\penalty0
  1908--1930, 2019.

\bibitem[Defazio(2019)]{defazio2019curved}
Aaron Defazio.
\newblock On the curved geometry of accelerated optimization.
\newblock In \emph{Advances in Neural Information Processing Systems}, pages
  1764--1773, 2019.

\bibitem[Diakonikolas and Orecchia(2019)]{diakonikolas2019approximate}
Jelena Diakonikolas and Lorenzo Orecchia.
\newblock The approximate duality gap technique: A unified theory of
  first-order methods.
\newblock \emph{SIAM Journal on Optimization}, 29\penalty0 (1):\penalty0
  660--689, 2019.

\bibitem[Drusvyatskiy(2017)]{drusvyatskiy2017proximal}
Dmitriy Drusvyatskiy.
\newblock The proximal point method revisited.
\newblock \emph{arXiv preprint: 1712.06038}, 2017.

\bibitem[Gasnikov and Nesterov(2018)]{gasnikov2018universal}
Alexander~Vladimirovich Gasnikov and Yu~E Nesterov.
\newblock Universal method for stochastic composite optimization problems.
\newblock \emph{Computational Mathematics and Mathematical Physics},
  58\penalty0 (1):\penalty0 48--64, 2018.

\bibitem[G{\"u}ler(1991)]{guler1991convergence}
Osman G{\"u}ler.
\newblock On the convergence of the proximal point algorithm for convex
  minimization.
\newblock \emph{SIAM Journal on Control and Optimization}, 29\penalty0
  (2):\penalty0 403--419, 1991.

\bibitem[Krichene et~al.(2015)Krichene, Bayen, and
  Bartlett]{krichene2015accelerated}
Walid Krichene, Alexandre Bayen, and Peter~L Bartlett.
\newblock Accelerated mirror descent in continuous and discrete time.
\newblock In \emph{Advances in Neural Information Processing Systems}, pages
  2845--2853, 2015.

\bibitem[Lan and Zhou(2018)]{lan2018optimal}
Guanghui Lan and Yi~Zhou.
\newblock An optimal randomized incremental gradient method.
\newblock \emph{Mathematical programming}, 171\penalty0 (1-2):\penalty0
  167--215, 2018.

\bibitem[Lee et~al.(2013)Lee, Rao, and Srivastava]{lee2013new}
Yin~Tat Lee, Satish Rao, and Nikhil Srivastava.
\newblock A new approach to computing maximum flows using electrical flows.
\newblock In \emph{Proceedings of ACM STOC}, pages 755--764, 2013.

\bibitem[Lessard et~al.(2016)Lessard, Recht, and Packard]{lessard2016analysis}
Laurent Lessard, Benjamin Recht, and Andrew Packard.
\newblock Analysis and design of optimization algorithms via integral quadratic
  constraints.
\newblock \emph{SIAM Journal on Optimization}, 26\penalty0 (1):\penalty0
  57--95, 2016.

\bibitem[Lewis and Wright(2016)]{lewis2016proximal}
Adrian~S Lewis and Stephen~J Wright.
\newblock A proximal method for composite minimization.
\newblock \emph{Mathematical Programming}, 158\penalty0 (1-2):\penalty0
  501--546, 2016.

\bibitem[Lin et~al.(2015)Lin, Mairal, and Harchaoui]{lin2015universal}
Hongzhou Lin, Julien Mairal, and Zaid Harchaoui.
\newblock A universal catalyst for first-order optimization.
\newblock In \emph{Advances in Neural Information Processing Systems}, pages
  3384--3392, 2015.

\bibitem[Martinet(1970)]{martinet1970regularisation}
Bernard Martinet.
\newblock R{\'e}gularisation d’in{\'e}quations variationnelles par
  approximations successives. rev. fran{\c{c}}aise informat.
\newblock \emph{Recherche Op{\'e}rationnelle}, 4:\penalty0 154--158, 1970.

\bibitem[Mokhtari et~al.(2019)Mokhtari, Ozdaglar, and
  Pattathil]{mokhtari2019unified}
Aryan Mokhtari, Asuman Ozdaglar, and Sarath Pattathil.
\newblock A unified analysis of extra-gradient and optimistic gradient methods
  for saddle point problems: proximal point approach.
\newblock \emph{arXiv preprint: 1901.08511}, 2019.

\bibitem[Moreau(1965)]{moreau1965proximite}
Jean-Jacques Moreau.
\newblock Proximit{\'e} et dualit{\'e} dans un espace {Hilbertien}.
\newblock \emph{Bulletin de la Soci{\'e}t{\'e} math{\'e}matique de France},
  93:\penalty0 273--299, 1965.

\bibitem[Nesterov(2007)]{nesterov2007modified}
Yu~Nesterov.
\newblock Modified {Gauss}-{Newton} scheme with worst case guarantees for
  global performance.
\newblock \emph{Optimisation methods and software}, 22\penalty0 (3):\penalty0
  469--483, 2007.

\bibitem[Nesterov(1983)]{nesterov1983method}
Yurii Nesterov.
\newblock A method for unconstrained convex minimization problem with the rate
  of convergence {$O(1/k^{2})$}.
\newblock In \emph{Doklady AN USSR}, volume 269, pages 543--547, 1983.

\bibitem[Nesterov(2018)]{nesterov2018lectures}
Yurii Nesterov.
\newblock \emph{Lectures on convex optimization}, volume 137.
\newblock Springer, 2018.

\bibitem[Rockafellar(1976)]{rockafellar1976monotone}
R~Tyrrell Rockafellar.
\newblock Monotone operators and the proximal point algorithm.
\newblock \emph{SIAM journal on control and optimization}, 14\penalty0
  (5):\penalty0 877--898, 1976.

\bibitem[Su et~al.(2016)Su, Boyd, and Candes]{su2014differential}
Weijie Su, Stephen Boyd, and Emmanuel~J Candes.
\newblock A differential equation for modeling {Nesterov's} accelerated
  gradient method: theory and insights.
\newblock \emph{JMLR}, 17\penalty0 (1):\penalty0 5312--5354, 2016.

\bibitem[Sutskever et~al.(2013)Sutskever, Martens, Dahl, and
  Hinton]{sutskever2013importance}
Ilya Sutskever, James Martens, George Dahl, and Geoffrey Hinton.
\newblock On the importance of initialization and momentum in deep learning.
\newblock In \emph{ICML}, pages 1139--1147, 2013.

\bibitem[Teboulle(2018)]{teboulle2018simplified}
Marc Teboulle.
\newblock A simplified view of first order methods for optimization.
\newblock \emph{Mathematical Programming}, 170\penalty0 (1):\penalty0 67--96,
  2018.

\bibitem[Tseng(2008)]{tseng2008accelerated}
Paul Tseng.
\newblock On accelerated proximal gradient methods for convex-concave
  optimization.
\newblock \emph{submitted to SIAM Journal on Optimization}, 2008.

\bibitem[Wang and Abernethy(2018)]{wang2018acceleration}
Jun-Kun Wang and Jacob~D Abernethy.
\newblock Acceleration through optimistic no-regret dynamics.
\newblock \emph{Advances in Neural Information Processing Systems}, 31, 2018.

\bibitem[Wibisono et~al.(2016)Wibisono, Wilson, and
  Jordan]{wibisono2016variational}
Andre Wibisono, Ashia~C Wilson, and Michael~I Jordan.
\newblock A variational perspective on accelerated methods in optimization.
\newblock \emph{PNAS}, 113\penalty0 (47):\penalty0 E7351--E7358, 2016.

\bibitem[Wilson et~al.(2016)Wilson, Recht, and Jordan]{wilson2016lyapunov}
Ashia~C Wilson, Benjamin Recht, and Michael~I Jordan.
\newblock A {Lyapunov} analysis of momentum methods in optimization.
\newblock \emph{arXiv preprint: 1611.02635}, 2016.

\end{thebibliography}

\appendix

\section{Deferred derivations}

\subsection{Deferred derivations from \S\ref{subsec:conv}} \label{app:agm}
Let us first derive \eqref{ineq:1:agm}. Applying Proposition~\ref{prop:per} with $\phi(x)=\ee_{t+1}[f(y_{t})+\inp{\nabla f(y_{t})}{x-y_{t}}]$ to \eqref{ppm:a}, we obtain:
\begin{align}
\phi(x_{t+1})-\phi(\xs) +\frac{1}{2}\norm{\xs-x_{t+1}}^2 +\frac{1}{2}\norm{ x_{t+1}-x_t}^2  -\frac{1}{2}\norm{\xs-x_t}^2\leq 0\,.  \label{ineq:basic}
\end{align}
Now from the convexity of $f$, it holds that $\phi(\xs) \leq \ee_{t+1} f(\xs)$. This together with the $L$-smoothness of $f$, it follows that \begin{align*}
    \phi(x_{t+1}) &= \ee_{t+1}[f(y_{t})+\inp{\nabla f(y_{t})}{z_{t+1}-y_{t}}+\inp{\nabla f(y_{t})}{x_{t+1}-z_{t+1}}] \\
    &\geq \ee_{t+1} \left[f(z_{t+1}) -\frac{L }{2}\norm{z_{t+1}-y_t}^2+\inp{\nabla f(y_{t})}{x_{t+1}-z_{t+1}} \right]\,.
\end{align*}  
Plugging these inequalities back to \eqref{ineq:basic} and rearranging, we obtain the following inequality:
\begin{align}
&\ee_{t+1}[f(z_{t+1})-f(\xs)] +\frac{1}{2}\norm{\xs-x_{t+1}}^2  -\frac{1}{2}\norm{\xs-x_t}^2 \nonumber\\
&\quad \leq-\frac{1}{2}\norm{ x_{t+1}-x_t}^2  +\ee_{t+1}\left[ \frac{L}{2}\norm{z_{t+1}-y_t}^2 +  \inp{\nabla f(y_{t})}{z_{t+1}-x_{t+1}}\right]\,. \label{exp:1} 
\end{align}
Now decomposing the inner product term  in \eqref{exp:1} into 
\begin{center}
$\ee_{t+1}\inp{\nabla f(y_{t})}{ z_{t+1}-y_t} +\ee_{t+1}\inp{\nabla f(y_{t})}{ y_t-x_t}+\ee_{t+1}\inp{\nabla f(y_{t})}{ x_t-x_{t+1}}$,    
\end{center} 
and using $x_{t+1}-x_t=-\ee_{t+1} \nabla f(y_t)$ and $z_{t+1}-y_t=-\nicefrac{1}{L} \nabla f(y_t)$ (which are \eqref{agm:b} and \eqref{agm:c}, respectively),  \eqref{exp:1} becomes $\left(\frac{\ee_{t+1}^2}{2}-\frac{\ee_{t+1}}{2L}\right) \norm{\nabla f(y_{t})}^2 + \ee_{t+1}\inp{\nabla f(y_{t})}{ y_{t}-x_t}$.  
Now, using the relation $y_{t}-x_t  =L\ee_t (z_t-y_t)$ (which is  \eqref{agm:a}), we obtain $\ffo$.
Thus, \eqref{ineq:1:agm} follows.

Next, \eqref{ineq:2:agm} readily follows from the $L$-smoothness and the convexity of $f$:
    \begin{align*}  
        f(z_{t+1}) -f(z_t) &= f(z_{t+1}) -f(y_{t})+f(y_{t})-f(z_t)\\ &\leq  \inp{\nabla f(y_{t})}{z_{t+1}-y_t}  +\frac{L}{2} \norm{z_{t+1}-y_t}^2+ \inp{\nabla f(y_{t})}{y_{t}-z_t}\\
        &\overset{(a)}{=}   -\frac{1}{2L} \norm{\nabla f(y_{t})}^2+ \inp{\nabla f(y_{t})}{y_{t}-z_t}=\fft,
    \end{align*}
    where ($a$) is due to $z_{t+1}-y_t=-\nicefrac{1}{L} \nabla f(y_t)$.

\subsection{Deferred derivations from \S\ref{sec:momentum}} \label{app:rmk}
We first derive \eqref{ineq:1:mom}.
By the updates \eqref{sim2}, we have $x_{t+1} = x_t -(\ee_t +\frac{1}{L})\nabla f(y_t)$. 
Letting $\eep_{t+1} :=\ee_t+\frac{1}{L}$, this relation can be equivalently written as:
\begin{align}
    \label{sim2:ppm} \textstyle x_{t+1} \leftarrow \argmin_{x} \left\{ f(y_t)+\inp{\nabla f(y_t)}{x-y_t}  + \frac{1}{2\eep_{t+1}} \norm{x-x_t}^2 \right\}
\end{align}
The rest is similar to \S\ref{app:agm}: we apply Proposition~\ref{prop:per} with $\phi(x)=\eep_{t+1}[f(y_{t})+\inp{\nabla f(y_{t})}{x-y_{t}}]$:
\begin{align}
\phi(x_{t+1})-\phi(\xs) +\frac{1}{2}\norm{\xs-x_{t+1}}^2 +\frac{1}{2}\norm{ x_{t+1}-x_t}^2  -\frac{1}{2}\norm{\xs-x_t}^2\leq 0\,. \label{ineq:basic2}
\end{align}
Now from the convexity, we have $\phi(\xs) \leq \eep_{t+1} f(\xs)$, and from the $L$-smoothness, we have 
\begin{align*}
    \phi(x_{t+1})&= \eep_{t+1} [f(y_{t})+\inp{\nabla f(y_{t})}{z_{t+1}-y_{t}}+\inp{\nabla f(y_{t})}{x_{t+1}-z_{t+1}}]\\
    &\geq \eep_{t+1} \left[f(z_{t+1}) -\frac{L }{2}\norm{z_{t+1}-y_t}^2+\inp{\nabla f(y_{t})}{x_{t+1}-z_{t+1}} \right]\,. 
\end{align*} 
Plugging these inequalities back to \eqref{ineq:basic2} and rearranging, we obtain the following inequality:
\begin{align*}
&\eep_{t+1}[f(z_{t+1})-f(\xs)] +\frac{1}{2}\norm{\xs-x_{t+1}}^2  -\frac{1}{2}\norm{\xs-x_t}^2 \\
&\quad \leq-\frac{1}{2}\norm{ x_{t+1}-x_t}^2  +\eep_{t+1}\left[ \frac{L}{2}\norm{z_{t+1}-y_t}^2 +  \inp{\nabla f(y_{t})}{z_{t+1}-x_{t+1}}\right] \\
&\quad =\frac{1}{2}\left(  -(L\ee_t+1)^2   +L\eep_{t+1}\right)\cdot \norm{z_{t+1}-y_t}^2 +  \eep_{t+1}\cdot \inp{\nabla f(y_{t})}{z_{t+1}-x_{t+1}} =\hho\,, 
\end{align*}
where the last line follows since $\norm{x_{t+1}-x_t} = (L\ee_t+1)\cdot \norm{z_{t+1}-z_t}$ (see Figure~\ref{fig:sim2}).

Next we derive \eqref{ineq:1:mom}. From the $L$-smoothness and the convexity of $f$:
    \begin{align*}  
        f(z_{t+1}) -f(z_t) &= f(z_{t+1}) -f(y_{t})+f(y_{t})-f(z_t)\\ &\leq  \inp{\nabla f(y_{t})}{z_{t+1}-y_t}  +\frac{L}{2} \norm{z_{t+1}-y_t}^2+ \inp{\nabla f(y_{t})}{y_{t}-z_t}\\
        &=  \frac{L}{2} \norm{z_{t+1}-y_t}^2+ \inp{\nabla f(y_{t})}{z_{t+1}-z_t}=\hht\,.
    \end{align*}

\subsection{Deferred derviations from \S\ref{sec:composite}} \label{app:sim}
Let us first derive \eqref{ineq:1:sim}. From convexity, we have $\pss(\xs) \leq \ee_{t+1} \fsi(\xs)$, and from the $L$-smoothness, we have the following lower bound:
\begin{align*}
    \pss(x_{t+1})&= \ee_{t+1}[f(y_{t})+\inp{\nabla f(y_{t})}{z_{t+1}-y_{t}}+\inp{\nabla f(y_{t})}{x_{t+1}-z_{t+1}} +\Psi(x_{t+1})]\\
    &\geq \ee_{t+1} \left[\fsi (z_{t+1}) -\frac{L }{2}\norm{z_{t+1}-y_t}^2+\inp{\nabla f(y_{t})}{x_{t+1}-z_{t+1}}  + \Psi(x_{t+1})-\Psi(z_{t+1}) \right]\,. 
\end{align*} 
Plugging these back to  \eqref{gen:1}, and using  the bound $- \breg{h}{x_{t+1}}{x_t}\leq -\frac{1}{2}\norm{x_{t+1}-x_t}^2$, \eqref{ineq:1:sim} follows.

Next, to derive \eqref{ineq:2:sim}, we use $L$-smoothness and the convexity of $f$ to obtain the following:
  \begin{align*}  
        \fsi(z_{t+1}) -\fsi(z_t) &\leq  f(z_{t+1})-f(y_t)+f(y_t)-f(z_t) +\Psi(z_{t+1})-\Psi(z_t)\\
        &\leq   \frac{L}{2} \norm{z_{t+1}-y_t}^2 +  \inp{\nabla f(y_{t})}{z_{t+1}-z_t} +\Psi(z_{t+1})-\Psi(z_{t})\,,
    \end{align*}
    which is precisely equal to $\ggt$.

\end{document}